\documentclass[11pt]{amsart}

\usepackage{amscd}
\usepackage{amssymb}
\usepackage{boxedminipage}
\usepackage{enumerate}
\usepackage{euscript}
\usepackage[all]{xy}

\theoremstyle{plain} \newtheorem{theorem}{Theorem}[section]
\theoremstyle{plain} \newtheorem{lemma}[theorem]{Lemma}
\theoremstyle{plain} \newtheorem{proposition}[theorem]{Proposition}

\newtheorem{corollary}[theorem]{Corollary}

\newtheorem{conjecture*}[]{Conjecture}

\newcommand{\nr}{\refstepcounter{theorem}  
                   \noindent {\thetheorem .}}
\newcommand{\defi}{\medskip \noindent {\it Definition \nr} }
\newcommand{\defifin}{\medskip}

\newcommand{\rem}{\medskip \noindent {\it Remark \nr} }
\newcommand{\remfin}{\medskip}

\newcommand{\gF}{\mathcal{F}}

\newcommand{\gE}{\mathcal{E}}

\newcommand{\gP}{\mathcal{P}}

\newcommand{\gL}{\mathcal{L}}

\newcommand{\gC}{\mathcal{C}}

\newcommand{\gG}{\mathcal{G}}
\newcommand{\gQ}{\mathcal{Q}}

\newcommand{\bD}{{\bf D}}
\newcommand{\bA}{{\bf A}}

\newcommand{\EH}{{\mathbf H}}

\newcommand{\Hom}{\text{Hom}}

\newcommand{\supp}{\text{supp}\,}

\newcommand{\sus}{\subseteq}

\newcommand{\pil}{\rightarrow}

\newcommand{\inpil}{\hookrightarrow}

\newcommand{\mto}[1]{\stackrel{#1}\longrightarrow}

\newcommand{\te}{\otimes}

\newcommand{\rk}{\text{\it rk} \,}

\newcommand{\tH}{\tilde{H}}

\newcommand{\Del}{\Delta}
\newcommand{\dl}{\Delta}

\newcommand{\lk}{\text{lk}}
\newcommand{\lkd}{\lk_\Delta}

\newcommand{\del}{\partial}

\newcommand{\Ga}{\Gamma}
\newcommand{\kan}{\omega}
\newcommand{\nat}{{\bf N}}

\newcommand{\hP}{\hat{P}}
\newcommand{\lkpo}{\lk_{P}}

\newcommand{\bb}{{\bf b}}

\begin{document}

\title [Cohen-Macaulay cell complexes ]
{Cohen-Macaulay cell complexes}
\author { Gunnar Fl{\o}ystad}
\address{ Matematisk Institutt\\
          Johs. Brunsgt. 12 \\
          5008 Bergen \\
          Norway}   
        
\email{ gunnar@mi.uib.no}

\begin{abstract}
We show that a finite regular cell complex with the intersection property
is a Cohen-Macaulay space iff the top enriched cohomology module
is the only nonvanishing one. 
We prove a comprehensive generalization of Balinski's
theorem on convex polytopes. 

An algebraic relation of K.Yanagawa implies that the (algebraic) enriched chain
and cochain complexes fit into a natural sixtuple of complexes
which in the simplicial case includes the resolution of the Stanley-Reisner
ring. A consequence turns out to be that there is no single generalization
of the Stanley-Reisner ring to cell complexes, but for Cohen-Macaulay 
cell complexes there is a generalization of the canonical module of the
Stanley-Reisner ring.
\end{abstract}

\maketitle


\section*{Introduction}

In \cite{Fl} we introduced enriched homology and cohomology
modules for simplicial complexes. They are modules 
over the polynomial ring $S = k[x_1, \ldots, x_n]$ where 
$\{1,2, \ldots,n\}$ are the vertices of the simplicial complex, and
their $S$-module ranks are equal to the $k$-vector space dimensions of
the corresponding reduced homology and cohomology groups.
We showed that a simplicial complex is Cohen-Macaulay iff there is  only one
nonvanishing enriched cohomology module, the top one.

In this paper we extended this to a class of cell complexes which
behaves nicely from a combinatorial viewpoint and includes polyhedral 
complexes. This is the class of finite regular cell complexes with
the intersection property. 
Enriched homology and cohomology modules may be defined equally well
for this class.
For simplicial complexes, being Cohen-Macaulay
is a topological property. We show that for a cell complex in the class
above, there is only one nonvanishing enriched cohomology  module
iff it is Cohen-Macaulay as a topological space.

In the class of Cohen-Macaulay simplicial complexes the classes 
of $l$-Cohen-Macaulay simplicial complexes constitute successively more
restricted classes. They correspond geometrically to successively 
higher connectivity of the simplicial complex. 
We show that this notion of being $l$-Cohen-Macaulay works well also for
the class of cell complexes we consider. We show that such a cell complex 
is $l$-Cohen-Macaulay iff its top enriched cohomology module can occur
as an $l-1$'t syzygy module in an $S$-free resolution. 
Also, its codimension $r$ skeleton will be 
$l+r$-Cohen-Macaulay. This generalizes results of \cite{Fl} for simplicial
complexes. It is also a comprehensive generalization of Balinski's theorem
for convex polytopes which says that the $1$-skeleton of a convex polytope
of dimension $d$ is $d$-connected. (Note that for a graph, being
$d$-connected is the same as being $d$-Cohen-Macaulay.)

\medskip
In \cite{Ya}, K.Yanagawa inroduced the notion of square free modules
over the polynomial ring $S$. Square free modules provides a natural 
setting for doing Stanley-Reisner theory. In \cite{Ya2} he defines two
dualities ${\bD}$ and ${\bA}$ on the category of bounded (algebraic)
complexes of free square free modules. He then shows the relation
\[  \bD \circ \bA \circ \bD \circ \bA \circ \bD \circ \bA (\gP) 
= \gP[-n] \]
where $\gP$ is a cochain complex and $\gP[-n]$ is the complex shifted $n$ steps
to the right. Starting from the enriched chain complex $\gE[-1]$ and 
applying ${\bD}$ and  ${\bA}$ we obtain a hexagon of complexes

\vskip 1mm
\begin{center}
\hskip -10mm
\xymatrix{ &  \gE[-1] \ar[dr]^{\bA}  & \\
  \gE^\vee[-1] \ar[ur]^{\bD} &    & \gG^\vee \ar[d]^{\bD} \\
\gF \ar[u]^{\bA[-n]} &    & \gG \ar[dl]^{\bA} \\
&    \gF^\vee \ar[ul]^{\bD} &.  \\ } 
\end{center}
\vskip 1mm

We show that when $\Gamma$ is simplicial the complex $\gF^\vee$ 
is the resolution of the Stanley-Reisner ring $k[\Gamma]$. If $\Gamma$
is not simplicial then $\gF^\vee$ will have more than one nonzero cohomology
module (which we describe). Thus one might say that although there is no
single generalization of the Stanley-Reisner ring to a cell complex,
there is a generalization of its resolution.

In case $\Gamma$ is Cohen-Macaulay and simplicial, the complex $\gF$ 
has only one nonvanishing homology module, the canonical module of the 
Stanley-Reisner ring $k[\Gamma]$. We show that for a Cohen-Macaulay cell
complex then $\gF$ still has only one homology module, which may then
be considered as a generalization of the canonical module.


\medskip

The organization of the paper is as follows. Section 1 consists of 
preliminaries on cell complexes, enriched homology and cohomology modules,
simplicial complexes and posets. Section 2 contains the main results, and
Section 3 the main bulk of the proofs of the results in Section 2. They are
formulated in the setting of posets. In Section 4 we describe 
the properties of the complexes in the hexagon above.

\medskip

{\bf Acknowledgement.} I thank the Mittag-Leffler Institute 
where I stayed part of the spring term 2005
and where the last part of this work was finished.

\section{Preliminaries}

\subsection{Cell complexes}

A {\it finite regular cell complex} $\Ga$ is a Haussdorff topological space 
$X$ together with a finite set of  closed subsets in $X$, 
homeomorphic to balls and called {\it cells} or {\it faces}, 
such that

\begin{itemize}
\item[i.] The interiors of the cells partition $X$.
\item[ii.] The boundary of each cell in is a union of
other cells.
\item[iii.] $\emptyset$ is a cell.
\end{itemize}

The dimension of a
cell is its topological dimension and the dimension of $\Ga$ is the
maximum of the dimension of the cells. We write $X = |\Ga|$.
The zero-dimensional cells are the vertices of $\Ga$ and we denote by $V$ 
the set of these. 

$\Ga$ is said to have {\it the intersection property} if
the intersection of any two faces is also a face. This means that the
face poset is a meet-semilattice. Two references for such cell complexes
are \cite{BB} and \cite{BK}.

\medskip

\noindent {\bf Note.} In this paper we shall only consider
finite regular cell complexes with the intersection property. To
avoid repeating this long phrase over and over, we shall simply use
the term {\it cell complex} as a short for {\it finite regular cell complex
with the intersection property}.

\medskip

If $R$ is a subset of the vertices, we denote by $\Ga_R$ the sub-complex
of $\Ga$ consisting of all faces whose vertex set is contained in $R$.
Also $\Ga_{V \backslash R}$ is denoted $\Ga_{-R}$.
We denote by $r$ the cardinality of $R$.

\subsection{Enriched homology and cohomology modules} \label{PreSecEnr}
The oriented augmented chain complex of $\tilde{\gC}(\Ga;k)$ of $\Ga$
consists of $\tilde{\gC}^i(\Ga;k)$ equal to the vector space over the field
$k$ spanned by
the $i$-dimensional faces, $\oplus_{\dim f = i} kf$. The differential
is given by
\[ f \mapsto \sum_{\dim f^\prime = i-1} \epsilon(f^\prime, f) f^\prime \]
where $\epsilon : \Ga \times \Ga \pil \{-1,0,1\} $ is a suitable incidence
function. The homology groups of this complex are the 
reduced homology groups $\tH_i(\Ga;k)$ of $\Ga$.  

\medskip
Now let $S$ be the polynomial ring $k[(x_v)_{v \in V}]$ in variables
indexed by the vertices $V$.
For a face $f$ let $m_f$ be the product of the variables indexed by the
vertices of $f$.
We now define the enriched homology modules of $\Ga$ in the same way
as we did for simplicial complexes in \cite{Fl}. Attaching the variable
$x_v$ to $v$, we may form the cellular complex $\gE(\Ga;k)$, see \cite{BS}, 
where $\gE_i(\Ga;k)$ is $\oplus_{\dim f = i} Sf$
The differential is given by 
\[ f \mapsto \sum_{\dim f^\prime = i-1} \epsilon(f^\prime, f) 
\frac{m_f}{m_{f^\prime}} f^\prime. \]
The {\it enriched homology module} $\EH_i(\Ga;k)$ (or just $\EH_i(\Ga)$)
is the $i$'th homology module of this complex. 
It is a module graded by $\nat^V$.
For $\bb$ in $\nat^V$ and $R$ the support of $\bb$, i.e. the set 
of nonzero coordinates, the graded part $\EH_i(\Ga)_{\bb}$ is the 
reduced homology $\tH_i(\Ga_R)$, see \cite{BS}. It follows as in 
\cite{Fl} that the $S$-module rank of $\EH_i(\Ga)$ is the dimension
of $\tH_i(\Ga)$ as a vector space over $k$.

\medskip
Let $\omega_S$ be the canonical module which is $S(-{\bf 1})$, i.e. the free
$S$-module with a generator in
multidegree ${\bf 1} = (1,1, \ldots, 1)$.
We define the {\it enriched cohomology module} $\EH^i(\Ga;k)$ (or just
$\EH^i(\Ga)$) as the $i$'th cohomology module of the dualized complex
\[ \gE(\Ga;k)^\vee = \Hom_S(\gE(\Ga;k), \omega_S). \] 

\subsection{ Simplicial complexes} \label{PreSusSim}
A simplicial complex $\dl$ on $[n] = \{ 1,2, \ldots, n\}$ is a family of 
subsets of $[n]$ such that if $X$ is in $\dl$ and $Y \sus X$ then $Y$ is 
in $\dl$. A standard reference is \cite{St}. Via its topological realization
it becomes a cell complex.

  For a subset $R$ of $[n]$, let $\dl_R$ denote the restricted simplicial
complex, consisting of faces that are subsets of $R$ and denote 
$\dl_{[n]\backslash R}$ by $\dl_{-R}$. Also let the link, $\lkd R$, be the
simplicial complex on $[n]\backslash R$ consisting of faces $Y$ such that
$Y \cup R$ is a face of $\dl$.

\medskip
 We may form the Stanley-Reisner ring $k[\dl]$ which is the polynomial
ring $k[x_1, \ldots, x_n]$ divided by the ideal generated by the
square free monomials corresponding to the non-faces of $\dl$.
We say that $\dl$ is Cohen-Macaulay (CM) if $k[\dl]$ is Cohen-Macaulay.
This is equivalent to the following homological criterion given by
Hochster \cite{Ho}
\begin{equation} \label{PreLigHoc}
 \tH_p(\dl_{-R}) = 0, \mbox{ when } p+r < \dim \dl.
\end{equation}
In \cite{Ba}, K.Baclawski defined a simplicial complex to be $l$-Cohen-Macaulay
if $\dl_{-R}$ is Cohen-Macaulay of the same dimension as $\dl$ for each
subset $R$ of cardinality less than $l$. By Hochsters criterion
(\ref{PreLigHoc}) this is equivalent to
\begin{equation}
\tH_p(\dl_{-R}) = 0, \mbox{ when } p+r < \dim \dl + l-1 \mbox{ and }
p < \dim \dl. \label{PreLigHlCM}
\end{equation}
If $k[\dl]$ is a Gorenstein ring and $\dl$ is not a cone, we say that
$\dl$ is Gorenstein*. This is equivalent to $\dl$ being $2$-Cohen-Macaulay and
$\tH^{\dim \dl}(\dl) = k$.

\subsection{Posets}
A standard reference for posets is \cite{StE}. Given a poset $P$ with
order relation $\leq$, an order ideal $J$ is a subset of $P$ such that
if $x$ is in $J$ and $y \leq x$, then $y$ is in $J$. 
A filter is subset $F$ such
that if $y$ is in $F$ and $y \leq x$, then $x$ is in $F$.
If  $R$ is subset of $P$, denote by $F(R)$ the filter generated by $R$,
consisting of $x$ in $P$ such that $x \geq r$ for some $r$ in $R$, and
for an element $x$ in $P$ denote by $P_{< x}$ the ideal of elements
less than $x$. 
An open interval $(x,y)$ consists of all $z$ strictly between $x$ and $y$,
$x < z < y$. 
We denote by $\hat{P}$ the poset $P$ with a bottom element $0$ and
top element $1$ adjoined.

The poset $P$ is graded if every maximal chain in $P$ has the same length
(the number of elements in the chain minus one). 
The length of such a maximal chain
is the rank of $P$, $\rk P$. For an element $x$ in $P$ the maximal chains
descending down from $x$, have the same length, the rank of $x$, $\rk x$.
If $\rk x = 0$, $x$ is called an atom.

The poset $P$ is a lattice if each pair $x$ and $y$ have a supremum $x \vee y$
and an infimum $x \wedge y$.

\medskip

For a poset $P$ we may form the order complex $\dl(P)$, a simplicial complex
consisting of all the chains in $P$. All the terminology for simplicial 
complexes may then be transferred to $P$. 
For simplicity we shall write $\tH_p(P_{-R})$ for $\tH_p(\dl(P)_{-R})$.
When $\Ga$ is a cell complex the nonempty cells form the face poset
$P(\Ga)$ with cells ordered by inclusion.
Then $\hat{P}(\Ga)$ will be a lattice (since we are assuming the
intersection property). It is a well-known fact that $\dl(P(\Ga))$ is
the complete barycentric subdivision of $\Ga$ and hence
as a topological space is homeomorphic to $|\Ga|$.

\section{Cohen-Macaulay cell complexes}
\label{CohSek}

\subsection {Criterion for being Cohen-Macaulay}

For a simplicial complex $\dl$, the property of being Cohen-Macaulay
is a topological property, \cite{Mu}. 
In fact, letting $X$ be the topological realization
$|\dl|$, then $\dl$ is CM iff

\begin{itemize}
\item[i.]  $H_p(X,X-\{x\};k) = 0$ for $p < \dim X$.
\item[ii.] $\tH_p(X) = 0$ for $p < \dim X$.
\end{itemize}
Hence $\dl$ and the face poset $P(\dl)$ are Cohen-Macaulay at the same
time. 
Now define a cell complex $\Ga$ to be Cohen-Macaulay if its face poset
$P(\Ga)$ is Cohen-Macaulay. Since $\dl(P(\Ga))$ is the complete 
barycentric subdivision of $\Ga$, this is a topological property of 
$\Ga$.
Recall the criterion (\ref{PreLigHoc}) of Hochster. We shall show that
this applies equally well as a criterion for a cell complex to 
be Cohen-Macaulay. 

\begin{theorem} \label{CMTheHov}
A cell complex $\Ga$ is Cohen-Macaulay if and only if 
$\tH_p (\Ga_{-R})$ vanishes whenever $R$ is a subset of the vertices $V$ 
with $p+r < \dim \Ga$.
\end{theorem}

\begin{proof}
The boundary of each cell $c$ is a sphere and so $P(\Ga)_{<c}$ is 
Gorenstein*. The only if direction now follows by applying Proposition 
\ref{PosProCMHo} and noting that $P(\Ga_{-R})$ is $P(\Ga)_{-F(R)}$.
The if direction follows by Proposition \ref{PosProHoCM}.
\end{proof}

\begin{corollary}
$\Ga$ is Cohen-Macaulay iff
the enriched cohomology modules $\EH^p(\Ga)$
vanish for $p < \dim \Ga$, or equivalently $\gE(\Ga)^\vee$ is a resolution of 
$\EH^{\dim \Ga}(\Ga)$.
\end{corollary}

\begin {proof}
The theorem just proven says 
$\tH_{\dim \Ga -i}(\Ga_{-R})$ vanishes for $r < i$. This means that
the enriched homology $S$-modules $\EH_{\dim \Ga -i}
(\Ga)$ have codimension $\geq i$. It is then a standard exercise in algebra
that this is equivalent to the dual complex $\gE(\Ga)^\vee$ of 
$\gE(\Ga)$ having no cohomology in (cohomological) degrees $< \dim \Ga$ by 
using the following (see \cite[18.4]{Ei} and its proof).

\vskip 2mm
\noindent {\bf Fact.} For an $S$-module
$M$ of codimension $r$, $Ext_S^p(M,\omega_S)$ vanishes for $p < r$, it has
codimension $r$ if $p = r$, and codimension $\geq p$ for $p \geq r$.
\end{proof}

\subsection{$l$-CM cell complexes}

Recall the notion in \ref{PreSusSim} of a simplicial complex being
$l$-Cohen-Macaulay. This notion turns out to work well for cell complexes
also.

\defi A cell complex $\Ga$ is $l$-Cohen-Macaulay if $\Ga_{-R}$ is
Cohen-Macaulay of the same dimension as $\Ga$ for each subset $R$ of 
the vertices $V$ of cardinality $< l$.
\defifin

The following, stated in \cite{Fl} for simplicial complexes, now holds.

\begin{theorem}
$\Ga$ is an $l$-CM cell complex iff the top cohomology module
$\EH^{\dim \Ga}(\Ga)$ can occur as an $l-1$'th syzygy module in a free
resolution.
\end{theorem}

\begin{proof}
Given the criterion of Theorem \ref{CMTheHov} applied to the restrictions
$\Ga_{-R}$ were $r < l$, the proof in \cite{Fl} carries over.
\end{proof}

Balinski's theorem for convex polytopes, see \cite{Zi}
says that the $1$-skeleton
of a $d$-dimensional polytope is $d$-connected. The following
gives a comprehensive generalization of this since a polytope, being
a ball, is $1$-Cohen-Macaulay.

\begin{corollary} If $\Ga$ is an $l$-CM cell complex, its codimension $r$
skeleton is $l+r$-CM.
\end{corollary}

\begin{proof}
Let $\Ga_{\leq \dim \Ga -r}$ be the codimension $r$ skeleton.
Then $\gE(\Ga_{\leq \dim \Ga -r})^{\vee}$ is the truncation of $\gE(\Ga)^\vee$
in cohomological 
degrees $\leq \dim \Ga -r$ and so the top cohomology module will
be an $l+r-1$'th syzygy module.
\end{proof}

Now it is known that for a simplicial complex, being $2$-CM is a topological
property, \cite{Wa}. (But being $l$-CM for $l \geq 3$ is not.)
The following shows the same for 
cell complexes.

\begin{theorem} 
A cell complex $\Ga$ is $2$-CM iff $P(\Ga)$ is a $2$-CM poset.
\end{theorem}

\begin{proof}
This follows by Propositions \ref{PosProHoCM2} and \ref{PosProCM2Ho}.
\end{proof}

\subsection{Gorenstein* cell complexes}

Recall that a simplicial complex $\dl$ is Gorenstein* iff it is $2$-CM and
$\tH_{\dim \dl}(\dl)$ is equal to $k$, see for instance the argument
in \cite[Thm. 3.1]{Fl}. We then define a cell complex $\Ga$ to be 
Gorenstein* iff it is $2$-CM and has $\tH_{\dim \Ga}(\Ga)$ equal to 
$k$. By the above theorem this is a topological property.

\begin{theorem} \label{CMTheGor}
$\Ga$ is a Gorenstein* cell complex iff $\Ga$ is Cohen-Macaulay and the
top cohomology module is a rank one torsion free $S$-module. It naturally
identifies as the ideal in $S$ generated be monomials $m_V/m_f$
where $f$ ranges over the facets of $\Ga$.
\end{theorem}

\begin{proof} This is completely analog to the proof of Theorem 3.1 in 
\cite{Fl}.
\end{proof}

\rem If $\Ga$ is a simplicial complex, the the quotient of $S$ by this ideal
is the Stanley-Reisner ring of the Alexander dual simplicial complex of $\Ga$.
\remfin

\section{Results on posets}

This section contains the essences of the proofs of the theorems
given in the previous Section \ref{CohSek}. We formulate the results here
in terms of posets, and this will be applied
to the face poset of a cell complex. Some ingredients will be used
repeatedly in the proofs and we inform on these first.

If $\dl$ is a simplicial complex and $v$ a vertex in $\dl$, 
let $\text{st}_{\dl} \{v\}$ consist of the faces of $\dl$ containing 
$v$. Thus $\lkd \{v\}$ is the intersection of $\text{st}_{\dl} \{v\}$ and
$\dl_{-\{v\}}$. Mayer-Vietoris for the pair $st_{\dl} \{v\}$ and 
$\dl_{-\{v\}}$ gives the long exact sequence
\[ \cdots \pil \tH_p(\lkd \{v\}) \pil \tH_p(\dl_{-\{ v\}}) \pil
\tH_p(\dl) \pil \cdots .\]
When $P$ is a poset and $x$ is in $P$, the link $\lkpo \{x\}$ is the
join $(0,x) * (x,1)$
consisting of subsets of $P$ which are unions of subsets of these
two intervals.
 The K\"unneth formula, then gives an 
isomorphism
\[ \tH_p( \lkpo \{ x\}) = \oplus_{i+j=p-1} \tH_i((0,x)) \te \tH_j((x,1)). \]

A wellknown criterion for $P$ to be CM is that for every interval $(x,y)$
in $\hP$, $\tH_p((x,y))$ vanishes unless $p$ is equal to the rank
of the interval, $\rk y - \rk x -2$.
We shall also in the following use the assumption that $P$ is 
Gorenstein*. This is equivalent to the additional requirement that
$\tH_p((x,y)) = k$ whenever $p$ equals the rank of the interval.
In particular every interval $(x,y)$ where $\rk y = \rk x +2$ contains
two elements.

\begin{lemma}  \label{PosLemAcy}
Let $P$ be a Gorenstein* poset and $x$ an element. 
Then $P_{-F(x)}$ is acyclic, i.e. $\tH_p(P_{-F(x)})$ is zero for each $p$.
\end{lemma}

\begin{proof}
First note that if $\dl$ is a Gorenstein* simplicial complex and $x$ is
a vertex in $\dl$, then $\tH_p(\dl_{-\{x\}})$ is zero for every $p$,
see for instance the argument in \cite[Thm.3.1]{Fl}.
Let $\{x \} \subsetneqq J^\prime \sus F(x)$ be an order ideal in $F(x)$.
Let $J^\prime = J \cup \{t\}$ where $t$ is a maximal element. We have
an exact sequence
\begin{equation} 
\tH_p(\lkpo\{t\}_{-J}) \pil \tH_p(P_{-J^\prime}) \pil \tH_p(P_{-J}) 
\label{PosLigSek1}
\end{equation}
where $\lkpo \{t\}_{-J}$ is the join of $(0,t)_{-F(x)}$ and $(t,1)$. By
induction on the sizes of $P$ and $J$ the outer groups of (\ref{PosLigSek1}) 
vanish and hence the middle one.
\end{proof}

\begin{lemma} \label{PosLemMax}
Let $P$ be Gorenstein* with $\hP$ a lattice, 
and let $R$ a subset of $P$ with at least 
two elements. Then there exists a maximal element
$s$ in $P$ such that $s \geq r$ for some $r$ in $R$, but not all $r$ in $R$.
\end{lemma}

\begin{proof} Let $t$ be the join
of all elements of $R$ and chose a minimal element $r$ in $R$. 
We construct successively larger elements $s_i$ which
fulfill the conditions above. First let $s_1 = r$. There are at least two
elements covering $r$. Not both of them are $\geq t$ since then their
meet $s_1 = r$ is $\geq t$. Let $s_2$ be one of them not $\geq t$. 
Continuing this process we finally arrive at a maximal $s$.
\end{proof}

\begin{proposition} \label{PosProCM2Ho}
 Assume $\hP$ is a graded lattice and $P_{<x}$ is
Gorenstein* for each $x$ in $P$. If $P$ is $2$-CM, then
$\tH_p(P_{-F(R)})$  is zero when $p+r \leq \rk P$ and $p < \rk P$.
\end{proposition}

\begin{proof}
We shall prove that $\tH_p(P_{-J^\prime})$ is zero when $R \subseteq J^\prime
\subseteq F(R)$ is an order ideal and $p+r \leq \rk P$ and $p < \rk P$.
We use induction on the sizes of $J^\prime$ and $R$. Note that the 
case $J^\prime = R$ is the hypothesis. So assume $J^\prime = J \cup \{t\}$
where $t$ is a maximal element in $J^\prime$ not in $R$. We have an
exact sequence
\begin{equation} \tH_{p+1}(P_{-J}) \mto{\eta} \tH_p(\lkpo \{t\}_{-J}) 
\pil \tH_p(P_{-J^\prime}) \pil \tH_p(P_{-J}). \label{PosLigSek2}
\end{equation}
By induction we may assume that the last term of this sequence is zero.
To achieve our goal we shall show that $\eta$ is surjective. 
We do this by induction on the sizes of $R$ and $J$.

Assume $R = \{x\}$ consists of one element. Then $\lkpo\{t\}_{-J}$ is the join
of $(0,t)_{-F(x)}$ and $(t,1)$. Using Lemma \ref{PosLemAcy} 
and the K\"unneth formula,
we get that the second term in (\ref{PosLigSek2}) 
is zero and so $\eta$ is surjective.

If $R^\prime = R \cap (0,t)$ contains less than $R$ elements, then
$\lkpo\{t\}_{-J}$ is the join of $(0,t)_{-F(R^\prime)}$ and $(t,1)$. 
If the second term in (\ref{PosLigSek2}) is nonzero then by the 
K\"unneth formula
$\tH_i((0,t)_{-F(R^\prime)})$ and $\tH_j((t,1))$ are nonzero for some pair
$i$ and $j$ with $i+j = p-1$. But then $j \geq \rk P - \rk t -1$ and by
inducion  $i+r^\prime \geq \rk t$. Thus $p+r^\prime \geq \rk P$.
This contradicts $p+r \leq \rk P$. Hence the second term in (\ref{PosLigSek2})
is zero and $\eta$ surjective.

So assume $r \geq 2$ and $R \subseteq (0,t)$. 
By Lemma \ref{PosLemMax}, 
there is a maximal $s$ in $(0,t)$ greater or equal to some
element in $R$ but not every element in $R$. Letting $J = J^{\prime \prime}
\cup \{s\}$ we get a diagram
\[ \begin{CD} \tH_{p+1}(\lkpo \{s\}_{-J^{\prime\prime}}) @>>>
\tH_{p+1}(P_{-J}) @>{\beta}>> \tH_{p+1}(P_{-J^{\prime\prime}}) \\
@VV{\eta^\prime}V @VV{\eta}V @VV{\eta^{\prime\prime}}V \\
\tH_p(\lkpo\{t,s\}_{-J^{\prime\prime}}) @>>> \tH_p(\lkpo \{t\}_{-J}) 
@>>> \tH_p(\lkpo \{t\}_{-J^{\prime\prime}}) \end{CD} \]
By induction $\eta^{\prime\prime}$ is surjective. We shall show that $\beta$
and $\eta^\prime$ are surjective. It then follows by a simple diagram
chase that $\eta$ is surjective. 
Let $R^\prime$ be $R \cap (0,s)$ which is strictly included in $R$. 
Now $\lkpo (\{s\}_{-J^{\prime\prime}})$ is the join of $(0,s)_{-F(R^\prime)}$
and $(s,1)$. Hence by induction and the K\"unneth formula it easily follows 
that $\tH_p(\lkpo \{s\}_{-J^{\prime\prime}})$ vanishes. Since this
continues the long exact sequence which $\beta$ belongs to, $\beta$ must
be surjective.

Now the long exact sequence to which $\eta^{\prime}$ belongs to is continued
by
\begin{equation} 
\tH_p(\lkpo \{s\}_{-J^{\prime\prime} \cup \{t\}}). \label{PosLigEfo}
\end{equation}
If nonzero, then by the K\"unneth formula
$\tH_i((0,s)_{-F(R^\prime)})$ and $\tH_j((s,1)_{-\{t\}})$ are nonzero for some
$i+j = p-1$. This gives $i+r^\prime \geq \rk s$ and, since $(s,1)$ is
$2$-CM since $P$ is (being $2$-CM is preserved by links), 
that $j+1 \geq \rk P - \rk s$. Then $p+r^\prime
\geq \rk P$, contradicting $p+r \leq \rk P$. Thus 
(\ref{PosLigEfo}) vanishes and so $\eta^\prime$
is surjective. This completes the proof.
\end{proof}

\begin{proposition} \label{PosProCMHo}
Assume $\hP$ is a graded lattice and $P_{<x}$ is Gorenstein* for every $x$
in $P$. If $P$ is Cohen-Macaulay, then 
$\tH_p(P_{-F(R)})$ is zero when $p+r < \rk P$.
\end{proposition}

\begin{proof}
We shall use induction on the size of an order ideal $R \subseteq J^\prime
\subseteq F(R)$ to prove that
$ \tH_p(P_{-J^\prime})$ is zero for $p+r < \rk P$. When $J=R$ this is the
hypothesis, so assume $J^\prime = J \cup \{t\}$ where $t$ is a maximal element
of $J^\prime \backslash R$. There is an exact sequence
\begin{equation}
\tH_p(\lkpo \{t\}_{-J}) \pil \tH_p(P_{-J^\prime}) \pil \tH_p(P_{-J})
\label{PosLigSek3}
\end{equation}
Letting $R^\prime$ be $R \cap (0,t)$, then $\lkpo \{t\}_{-J}$ is the join
of $(0,t)_{-F(R^\prime)}$ and $(t,1)$. By Proposition \ref{PosProCM2Ho}
we know that
$\tH_i((0,t)_{-F(R^\prime)})$ nonzero implies $i+r^\prime \geq \rk t$.
Also $\tH_j((t,1))$ nonzero implies $j \geq \rk P - \rk t -1$. Hence
the K\"unneth formula gives $\tH_p(\lkpo \{t\}_{-J})$ is nonzero only if 
$p+r^\prime = i+j+r^\prime +1 \geq \rk P$. Since we are assuming 
$p+r < \rk P$, the outer terms in (\ref{PosLigSek3}) vanish and so also
the middle term.
\end{proof}

\begin{proposition} \label{PosProHoCM}
Assume $\hP$ is a graded lattice with 
$P_{<x}$ Gorenstein* for all $x$ in $P$. If $\tH_p (P_{-F(R)})$ is zero
when $p+r < \rk P$ and  $R$ consists of atoms, 
then  $P$ is Cohen-Macaulay.
\end{proposition}

\begin{proof}
We first show that $\tH_p(P_{-F(R)})$ is zero for any $R$ with $p+r < \rk P$.
We shall use induction on the rank of a maximal element in $R$ and
the number of elements in $R$ with this maximal rank.
So let $x$ in $R$ be of maximal rank in $R$, supposed to be $> 0$.
There is thus in $\hP$ an element $w$ with $\rk_{\hP} x = \rk_{\hP} w +2$.
Since $P_{< x}$ is Gorenstein*, there are at least two elements strictly
between $w$ and $x$, say $y$ and $z$ with $y \wedge z = w$.
Let $R_0 = R\backslash \{x\}$. Then 
\[ F(R_0 \cup \{x\}) = F(R_0 \cup \{y\}) \cap F(R_0 \cap \{z\}).\]
Hence the order complex $\dl(P)_{-F(R_0 \cup \{x\})}$ is the union of
$\dl(P)_{-F(R_0 \cup \{y\})}$ and $\dl(P)_{-F(R_0 \cup \{z\})}$, and 
$\dl(P)_{-F(R_0 \cup \{y,z\})}$ is the intersection of these two complexes.
The Mayer-Vietoris sequence 
\[ \tH_p(P_{-F(R_0 \cup \{y\})}) \oplus \tH_p(P_{-F(R_0 \cup \{z\})})
\pil \tH_p(P_{-F(R)}) \pil \tH_{p-1}(P_{-F(R_0 \cup \{y,z\})}) \]
then gives by induction that the middle term vanishes for $p+r < \rk P$.

\medskip

In order to show that $P$ is CM, it will be enough to show that 
$\tH_j((y,1))$ vanishes for 
$j < \rk P - \rk y -1$. Since $\tH_i((0,y))$ is $k$ if $i = \rk y -1$ and
zero otherwise, this is by the K\"unneth formula the same as showing that
$\tH_p(\lkpo \{y\})$ vanishes for $p < \rk P -1$. From the exact sequence
\[ \tH_{p+1}(P) \pil \tH_p(\lkpo \{y\}) \pil \tH_p(P_{-\{y\}}) \pil \tH_p(P) \]
it follows that we can show that the third term is zero.

Let $\{y \} \subseteq J^\prime \sus F(y)$ be an order ideal. 
We show by descending
induction on the size of $J^\prime$ that $\tH_p(P_{-J^\prime})$ is zero.
Let $J^\prime = J \cup \{t\}$ where $t > y$ is a maximal element in $J^\prime$.
There is an exact sequence
\[ \tH_p(\lkpo \{t\}_{-J}) \pil \tH_p(P_{-J^\prime}) \pil \tH_p(P_{-J})
\pil \tH_{p-1}(\lkpo \{t\}_{-J}). \]
Here $\lkpo \{t\} $ is the join of $(0,t)_{-F(y)}$ 
which is acyclic by Lemma \ref{PosLemAcy},
and $(t,1)$. Hence $\tH(P_{-J})$ vanishes
by induction.
\end{proof}

\begin{proposition} \label{PosProHoCM2}
Assume $\hP$ is a graded lattice such that $P_{< x}$ is Gorenstein*
for each $x$ in $P$. If $\tH_p(P_{-F(R)})$ is zero when $R$ consists
of atoms, $p+r \leq \rk P$ and $p < \rk P$, then $P$ is 2-CM.
\end{proposition}

\begin{proof}
It follows exactly as in Proposition \ref{PosProHoCM} 
that $\tH_p(P_{-F(R)})$ is zero
for any $R$ with $p+r \leq \rk R$ and $p < \rk R$. 

We shall show that $\tH_p(P_{-R})$ is zero when $p+r \leq \rk P$ and
$p < \rk P$. By Hochsters criterion 
(\ref{PreLigHlCM})
this gives $P$ is $2$-CM.
 We may assume $R$ is nonempty.
 Let $R \subseteq J \subseteq F(R)$ be an order
ideal. We prove by descending induction on the size of $J$ that 
$\tH_p(P_{-J})$ is zero when $p+r \leq \rk P$ and $p < \rk P$.
Assume $J \subsetneqq F(R)$ and let $J^\prime = J \cup \{t\}$
be a strictly larger order ideal in $F(R)$. There is an exact sequence
\begin{equation}
\tH_p(P_{-J^\prime}) \pil \tH_p(P_{-J}) \pil \tH_{p-1}(\lkpo \{t\}_{-J})
\label{PosLigHo2CM}
\end{equation}
If the last term is nonzero, then $\tH_i((0,t)_{-F(R)})$ 
and $\tH_j((t,1))$ are nonzero for some $i,j$ with $i+j = p-2$. 
By induction $i+r \geq \rk t$ and by Proposition \ref{PosProHoCM}, 
$j \geq \rk P - \rk t -1$. 
Hence $p-2+r \geq \rk P -1$ or $p+r \geq \rk P + 1$, which is against
assumption. Hence the outer terms of (\ref{PosLigHo2CM}) vanish and
so the middle term.
\end{proof}

\section{The hexagon of complexes}

In this section we show that although there is no single generalization
of the Stanley-Reisner ring to cell complexes, the resolution of the
Stanley-Reisner ring, a complex of free $S$-modules, may be generalized
to a complex of free $S$-modules associated to any cell complex.
This complex of free $S$-modules together with
the enriched chain and cochain complexes fit into a natural sixtuple
of complexes associated to any cell complex and we shall describe
the properties of this sixtuple. A consequence is that the canonical
module of the Stanley-Reisner ring of a Cohen-Macaulay simplicial complex
may be generalized to a module associated to any Cohen-Macaulay 
cell complex.

\subsection{Square free modules and complexes}

In \cite {Ya}, Yanagawa introduced the notion of square free modules 
over the polynomial ring $S$. An $\nat^V$-graded module $M$ over $S$ is 
square free
if the multiplication map $M_{\bb} \mto{\cdot x_v} M_{\bb + e_v}$, where
$e_v$ is the $v$'th coordinate vector, is a bijection when $v$ is contained
in the support of $\bb$, i.e. the set of non-vanishing coordinates.

\medskip
If $R$ is a subset of $V$ we shall, where appropriate identify
$R$ with the multidegree $\bb$ which is the characteristic vector of $R$.
For a square free module $M$, independently T.R\"omer \cite{Ro} and
E.Miller \cite{Mi} defined its Alexander
dual $M^*$ as follows. 
For a subset $R$ of $V$, 
$(M^*)_R$ is the dual $\Hom_k(M_{R^c},k)$.
If $v$ is not in $R$ the multiplication
\[ (M^*)_R \mto{\cdot x_v} (M^*)_{R \cup \{v\}}  \]
is the dual of the multiplication
\[ M_{(R\cup \{v\})^c} \mto{\cdot x_v} M_{R^c}. \]
By obvious extension this defines $(M^*)_\bb$ for all $\bb$ in $\nat^V$
and all multiplications.

\medskip

A free square free module is a module
\[ \bigoplus_{R \sus V} S \te_k B_R \]
where $B_R$ is a vector space of multidegree $R$. In a cochain complex $\gP^\cdot$
of free square free modules with term
\[ \gP^i = \bigoplus_{R \sus V} S \te_k B^i_R \]
the spaces $B_R^i$ are called the {\it Betti spaces} of $\gP^\cdot$. The 
differentials in such a complex are to be homogeneous of degree zero and
the complex is called minimal if the differentials in $k \te_S \gP^\cdot$ vanish.
The complex $\gP^\cdot$ may be shifted $n$ steps to the left (resp. right
if $n < 0$) to a complex $\gP^\cdot[n]$ where
$\gP[n]^i$ is equal to $\gP^{i+n}$. 
We may also consider a cochain complex $\gP^\cdot$ as a chain complex
$\gP_\cdot$ by letting $\gP_i$ be $\gP^{-i}$. Note that $\gP^\cdot[n]$ then is 
$\gP_\cdot[-n]$. 

\medskip
On the category of complexes of free square free modules there are two
functors $\bD$ and $\bA$, studied in \cite{Ya2}. The functor $\bD$ is 
given by
\[  \bD(\gP^\cdot) = \Hom_S(\gP^\cdot, \omega_S). \]
The functor $\bA$ is given by first taking the Alexander dual complex
$\gP^{\cdot*}$. Unfortunately this is not a complex of {\it free} square free modules.
What
we do is to take a minimal resolution $\gQ^\cdot$ of this complex. I.e $\gQ^\cdot$
is a minimal complex of free modules together with a morphism 
$\gQ^\cdot \pil \gP^{\cdot*}$ which induces an isomorphism on cohomology. We then write
$\bA(\gP^\cdot) = \gQ^\cdot$.

When $\gP^\cdot$ is a minimal cochain complex, Yanagawa in \cite{Ya2} 
shows the fundamental 
isomorphism 
\begin{equation} \bD \circ \bA \circ \bD \circ \bA \circ \bD \circ \bA (\gP^\cdot) 
= \gP^\cdot[-n]. \label{HexLigDA}
\end{equation}

\rem In \cite{Ya2} $\bD$ is defined slightly diffferently. There 
$\bD$ is defined as \[ \bD(\gP^\cdot) = \Hom_S(\gP^\cdot, \omega_S)[n]\] 
which is natural in a general algebraic setting.
Then the relation (\ref{HexLigDA}) is
\[ \bD \circ \bA \circ \bD \circ \bA \circ \bD \circ \bA (\gP^\cdot) 
= \gP^\cdot[2n]. \]
\remfin

By starting with the enriched chain complex shifted one step to the left, 
$\gE[-1]$, we get by 
(\ref{HexLigDA}) a hexagon


\vskip 3mm
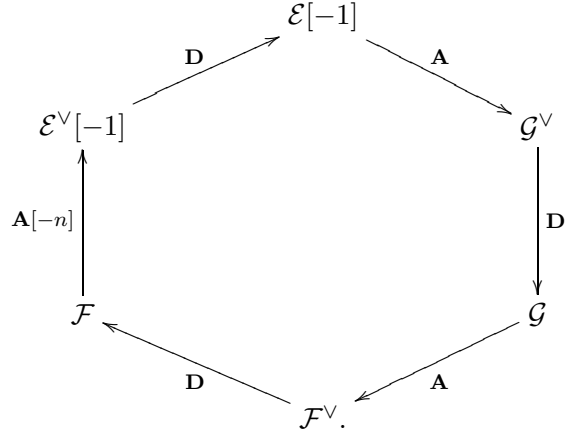
\begin{center}
\hskip -2mm
\xymatrix{ & & \gE[-1] \ar[drr]^{\bA}&  & \\
  \gE^\vee[-1] \ar[urr]^{\bD} & & &   & \gG^\vee \ar[dd]^{\bD} \\
 & & & & \\
\gF \ar[uu]^{\bA[-n]} & & &   & \gG \ar[dll]^{\bA} \\
&   & \gF^\vee. \ar[ull]^{\bD}& &  \\  }

\end{center}
\vskip 3mm


We consider $\gE[-1],\gF$, and  $\gG$ as chain complexes and their duals by 
$\bD$ as cochain complexes.
Now we proceed to study these complexes.

\subsection{Homology and Betti spaces.} For
a complex $\gP^\cdot$ of free $S$-modules define its $i$'th linear strand
$\gP^\cdot_{\langle i \rangle}$ to be given by 
\[ \gP^{j}_{\langle i \rangle} = \bigoplus_{|R| = i-j} S \te_k B^j_R. \]

For a square free module $M$, one may define a complex $\gL(M)$ (see
\cite[p.9]{Ya2} where it is denoted by $\gF(M)$) by 
\[\gL^i(M) = \bigoplus_{|R| = i} (M_R)^\circ \te_k S\] 
where $(M_R)^\circ$ is $M_R$ but considered to have multidegree
$R^c$. The 
differential is 
\[ m^\circ \te s \mapsto \sum_{j \not \in R} (-1)^{\alpha(j,R)} (x_jm)^\circ
\te x_js \]
where $\alpha(j,R)$ is the number of $i$ in $R$ such that $i < j$
after putting some total order on $V$.

The following is \cite[Thm. 3.8]{Ya2}.

\begin{proposition} \label{HexProDA} Let $\gQ^\cdot$ be $\bD \circ \bA (\gP^\cdot).$
Then 
\[ \gL(H^i(\gQ^\cdot))[n-i]\] is the $i$'th linear strand of $\gP^\cdot$.
\end{proposition}

\begin{corollary} \label{HexCorADA}
Let $\gQ^\cdot$ be $\bA \circ \bD \circ \bA (\gP^\cdot)$. Then 
\[ \Hom_S(\gL(H^{-i}(\gQ^\cdot)), \omega_S)[-i] \]
is the $i$'th linear strand in $\gP^\cdot$. (And similarly with $\gP^\cdot$ and
$\gQ^\cdot$ interchanged.)
\end{corollary}

More informally, for a pair $\gP^\cdot$ and $\gQ^\cdot$ on opposite corners of the
hexagon, the Betti spaces of $\gP^\cdot$ corresponds to the
homology spaces of  $\gQ^\cdot$ and vice versa.

\medskip
Let $k^i[\Ga]$ be the square free
$S$-module given by 
$k^i[\Ga]_F = kF$ if $F$ is the vertex set of some face $f$ 
where the cardinality of $F$ is $\dim f +i+1$ and $k^i[\Ga]_F = 0$ if 
$F$ is not so. For two faces
$f^\prime$ and $f$ with vertex sets $F\cup \{v\}$ and $F$ the multiplication
\[k^i[\Ga]_F \mto{\cdot x_v} k^i[\Ga]_{F\cup \{v\}} \]
is given by sending $F$ to $F \cup \{v\}$.
Note that when $\Ga$ is a simplicial complex, $k^0[\Ga]$ is the 
Stanley-Reisner ring and $k^i[\Ga]$ is zero for $i$ not $0$.

 The main observation in this section is the following which shows that
$\gF^\vee$ is a generalization of the resolution of the Stanley-Reisner
ring of a simplicial complex.

\begin{theorem} \label{HexThmF}
$H^{-i}(\gF^\vee) = k^i[\Gamma]$
\end{theorem}

\begin{proof} This follows by Corollary \ref{HexCorADA} together with the fact 
that the linear strands in $\gE[-1]$ are exactly
\[ \Hom_S(\gL(k^i[\Gamma]), \omega_S)[i]. \]
\end{proof}

\begin{corollary} When $\Gamma$ is simplicial then $\gF^\vee$ is 
the resolution of the 
Stanley-Reisner ring and $\gG$ is the resolution of the Stanley-Reisner ideal 
$I_{\Gamma^*}$ of the Alexander
dual simplicial complex $\Gamma^*$.
\end{corollary}

\begin{proof}
When $\Gamma$ is simplicial then $k^i[\Gamma]$ is zero for $i > 0$
and $k^0[\Gamma]$ is the Stanley-Reisner ring.  
Hence $\gF$ is a resolution of the Stanley-Reisner ring.
Since the Alexander dual
module of $k[\Gamma]$ is exactly $I_{\Gamma^*}$, the homology of $\gG$ is
this ideal in homological degree zero and vanishes elsewhere. So $\gG$ is
a resolution of this ideal.
\end{proof}

Now the Betti spaces of the enriched chain complex are of course given
by 
\[  B_{iF}(\gE[-1]) = kF \]
when $F$ is the vertices of a face of  $\Gamma$ of dimension $i-1$. 
The homology of $\gE[-1]$ is, as noted in Section \ref{PreSecEnr}, given by
\[ H_i(\gE[-1])_{\bb} = \tH_{i-1}(\Gamma_{\supp \bb}). \]
On the other hand the cohomology of the enriched cochain complex 
$\gE^\vee$ is by \cite{Fl} in the simplicial case given by
\[ H^i(\gE^\vee[-1]) = \tH^{i-1-|(\supp{\bb}^c)|}(\lk_\Gamma (\supp \bb)^c).\]
Using Proposition \ref{HexProDA} 
together with the fact that the homologies of 
$\gP^\cdot$ and $A(\gP^\cdot)$ are Alexander dual modules, and the Betti spaces of $\gP^\cdot$ 
and $D(\gP^\cdot)$ are Alexander dual spaces, we get in the simplicial case the
following table over the homology and Betti spaces.

\begin{tabular}{c||c|c}
&&\\
Simplicial complex & Homology/cohomology & Betti spaces \\
$\Del$ & $H_i(-)_{\bb} /H^i(-)_{\bb}$ & $B_{iF}(-)/B^i_F(-)$ \\
&&\\

\hline 
\hline
&&\\
Enriched chain complex & Restriction & Face \\
$\gE[-1]$ & $\tH_{i-1}(\Del_{\supp \bb})$ & $kF$ \\
&&\\

\hline
&&\\
  & Restriction & Link \\
$\gG^\vee $ & $\tH^{i-1} (\Del_{(\supp \bb)^c})$ & $\tH^{i-1}( \lk_\Del F)$ \\
&&\\

\hline
&&\\
Resolution of $I_{\Del^*}$ & Face ($i=0$) & Link \\
$\gG$ & $k(\supp \bb)^c$ & $\tH_{i-1}(\lk_\Del F^c)$ \\
&&\\

\hline
&&\\
Resolution of $k[\Del]$ & Face ($i=0$) & Restriction \\
$\gF^\vee$ & $(k\cdot\supp \bb)^*$ & $\tH^{i-1+|F|}(\Del_F)$ \\
&&\\

\hline
&&\\
  &  Link & Restriction \\
$\gF$ & $\tH_{n+i-1-|(\supp \bb)|}(\lk_\Del (\supp \bb))$ &
$\tH_{i-1+|F|}(\Del_{F^c})$ \\
&&\\

\hline
&&\\
Enriched cohain complex & Link & Face \\
$\gE^\vee$ & $\tH^{i-1-|(\supp \bb)^c|} (\lk_\Del (\supp \bb)^c)$ & $(kF^c)^*$ \\
&&
\end{tabular}

\subsection{Linear complexes.}
A theorem of Eagon and Reiner \cite{ER} says that when $\Gamma$ is
simplicial then the resolution $\gG$ of $I_{\Gamma^*}$ is linear iff
$\Gamma$ is Cohen-Macaulay. The following generalizes this.

\begin{theorem}
$\gG$ is linear iff $\Gamma$ is a Cohen-Macaulay cell complex.
\end{theorem}

\begin{proof}
The linear strands of $\gG$ corresponds to the cohomology of $\gE^\vee[-1]$.
But $\Gamma$ is a Cohen-Macaulay cell complex iff $\gE^\vee[-1]$ has only one
nonvanishing cohomology module and hence $\gG$ only one linear strand.
\end{proof}

We then readily get the following.
\begin{itemize}
\item[] $\gE$ and $\gE^\vee$ are linear iff $\Gamma$ is simplicial.
\item[] $\gG$ and $\gG^\vee$ are linear iff $\Gamma$ is Cohen-Macaulay.
\item[] $\gF$ and $\gF^\vee$ are linear iff $\Gamma$ is a simplex on the
vertices it contains.
\end{itemize}

\subsection{Generalizations of the canonical module.}

When $\Gamma$ is simplicial and Cohen-Macaulay the complex $\gF$ has only
one cohomology module, the canonical module $\omega_{k[\Gamma]}$, and
this is the Alexander dual module of the only nonvanishing cohomology module
of $\gE^\vee[-1]$, the top enriched cohomology module of $\Gamma$. 
For a Cohen-Macaulay cell complex we may therefore consider the Alexander
dual of its top cohomology module, let us denote it by
$\omega_\Gamma$, as a generalization of the canonical module.

\subsection{Convex polytopes.}

If $\Ga$ is a convex polytope with boundary $\del \Ga$ then
$\EH^{\dim \del \Ga}(\del \Ga)$ identifies as an ideal in $S$ by 
Theorem \ref{CMTheGor}.
Taking the Alexander dual of the inclusion 
$\EH^{\dim \del \Ga}(\del \Ga) \inpil S$ we get a surjection
$\kan_{\del \Ga} \leftarrow S$ and so $\kan_{\del\Ga}$ identifies as
a quotient ring of $S$. When $\Ga$ is simplicial, the canonical
module $\omega_{k[\del\Ga]}$ is isomorphic to the Stanley-Reisner ring
$k[\del\Ga]$.  Thus in general for convex polytopes $\Ga$ one might consider
$\kan_{\del \Ga}$ as a generalization of the Stanley-Reisner ring, although
I would not consider it as  
the fully natural viewpoint since one then should consider the
complex $\gF^\vee$ instead of $\gF$. Theorem \ref{HexThmF} also gives that the 
quotient ring 
$\kan_{\del\Ga}$ will be Cohen-Macaulay only when $\Ga$ is simplicial.

Concretely the ring $\kan_{\del\Ga}$ may be described as the square free ring
such that $(\kan_{\del\Ga})_F$ is $k$ when $F$ is contained in a facet of 
$\del \Ga$
and $0$ when $F$ is not contained in a facet. It is the quotient of $S$ by
the ideal generated by the
$m_F$ where $F$ is not contained in a facet.

\end{document}